\documentstyle[12pt]{article}
\newcommand{\const}{\mathop{\rm const}\limits}

\begin{document}

\begin{center}

         {\bf CRITERION FOR CONVERGENCE ALMOST EVERYWHERE, \\

\vspace{3mm}

          with applications} \par

\vspace{4mm}

 $ {\bf E.Ostrovsky^a, \ \ L.Sirota^b } $ \\

\vspace{4mm}

$ ^a $ Corresponding Author. Department of Mathematics and computer science, Bar-Ilan University, 84105, Ramat Gan, Israel.\\

E-mail:\ eugostrovsky@list.ru\\

\vspace{4mm}

$ ^b $  Department of Mathematics and computer science. Bar-Ilan University,
84105, Ramat Gan, Israel.\\

E-mail:\ sirota3@bezeqint.net

\vspace{3mm}
                    {\sc Abstract.}\\

 \end{center}

 \vspace{3mm}

  \ We derive the necessary and sufficient condition for almost sure convergence of the sequence of measurable functions,
and consider some applications in the theory of Fourier  series and in the theory of random fields. \par

\vspace{4mm}

\par

\vspace{4mm}

\section{Introduction, Notations. Statement of problem.}

 \vspace{4mm}

  \ Let $ (X = \{ x \},M, \mu) $ be measurable space equipped with sigma-finite non zero measure measure $  \mu, \hspace{5mm}
(B = \{b\}, || \ \cdot \ || =  || \ \cdot \ ||B) $ be a (complete) separable Banach space relative the norm function $  || \ \cdot \ ||
= || \ \cdot \ ||B, $  not necessary to be separable (or reflexive). \par
  \ Let also

  $$
   F = \left\{ \ f(x) = f_{\infty}(x), \ \{f_n(x) \}, \ n = 1,2,\ldots; \ x \in X  \right\}
  $$
be a numbered family of  measurable  functions from the set $  X $ into the space $  B: $

 $$
f_n: X \to B, \ n \in \{ \infty \} \cup \{1,2,  \ldots \}.
 $$

\vspace{3mm}

{\bf Our goal in this short article is finding of the necessary and sufficient condition on the family $  F  $
in the integrals terms for almost sure convergence of the sequence
$ f_n \to f_{\infty}  $ or simple  $ \exists \lim_{n \to \infty} f_n \ a.e.:  $  }

$$
\mu \{x, \ x \in X, \ \lim_{n \to \infty} f_n(x) \ne f(x) \} = 0. \eqno(1.0)
$$

\ {\bf Remark 1.1.} \ We imply in the equality (1.0) that if the limit $  \lim_{n \to \infty} f_n(x)  $
can not exists, but on the set  with zero measure.\par

\vspace{3mm}

 \ Immediate predecessor of general case besides the special cases: martingales, monotone sequences etc.
 is the preprint \cite{Ostrovsky301}. We intend to generalize the results obtained therein.\par

 \ There are many applications of solution of this problem in the theory of Fourier  series (and integrals)
 \cite{Lukomskii1}, \cite{Reyna1} and other orthogonal ones \cite{Alexits1},
 theory of Probability, \cite{Billingsley1}, \cite{Billingsley2}, in particular, in the theory of martingales \cite{Hall1},
statistics \cite{Pawlak1} etc.\par

 \vspace{4mm}

\section{ Main Result: Convergence of a  Sequence of a measurable Functions almost everywhere.}

 \vspace{3mm}

 \hspace{3mm} {\bf Lemma 2.1.} There exists a {\it probabilistic} measure $  \nu $ defined on at the same sigma-field $  M $
which is equivalent to the source measure $  \mu  $ in the Radon-Nikodym sense,  i.e. such that the following
implication there holds

$$
\forall A \in M \ \Rightarrow [ \mu(A) = 0 \Leftrightarrow \nu(A) = 0]. \eqno(2.1)
$$

 \vspace{3mm}

 \ {\bf Proof.}  The case of boundedness of the measure $  \mu(\cdot)  $ is trivial; suppose therefore $  \mu(X) = \infty. $ \par

 \ As long as the measure $  \mu  $ is sigma - finite, there exist a countable family of disjoint
measurable sets $  \{X_m \}, \ m = 1,2,3, \ldots, \ X_m \in M, \ l \ne m \Rightarrow  X_m \cap X_l = \emptyset, $  such that

$$
 0 < \mu(X_m) < \infty, \hspace{5mm}  \cup_{m=1}^{\infty} X_m = X, \eqno(2.2)
$$
so that the family $  \{X_m \}  $ forms the partition of whole set $  X. $\par

 \ Define for any set $  A \in M $

$$
\nu(A) \stackrel{def}{=} \sum_{m=1}^{\infty}  \frac{ \mu(A \cap X_m)}{2^m \ \mu(X_m)}. \eqno(2.3)
$$

 \ Obviously, $  \nu(\cdot) $ is sigma-additive probability  $  (\nu(X) = 1)  $ measure which is completely
equivalent to the initial measure $   \mu. $ \par

\vspace{3mm}

 \ As a consequence: the $ \ \mu \ - $ complete convergence of the sequence of measurable functions $ \{f_n(x) \} $ is
equal to the convergence with $  \ \nu \ - \ $ measure one.\par

\vspace{3mm}

 {\it  Therefore, we can and will assume without loss of generality that the initial measure $   \mu   $ is
 probabilistic, i.e. $  \mu(X) = 1. $ } \par

 \ Further, we can reduce our problem again without loss of generality  passing to the sequence $  g_n = f_n - f_{\infty}  $
to the case when $ f_{\infty}  = 0.$  \par

 \ A probability language: given the sequence of  random variables $  \xi_n, \ \xi, \\ \ n = 1,2,\ldots  $  with values in certain separable
Banach space $ B; $  find the necessary and sufficient condition for the following  convergence with probability one:

$$
{\bf P} ( \lim_{n \to \infty} || \ \xi_n - \xi \ ||B = 0 ) = 1. \eqno(2.4)
$$

 \ One can assume as before $  \xi = 0. $  \par

 \ Let as introduce the Banach space $  c_0(B), $ also separable, consisting on all the sequences  $  y = \{ y(n) \}  $
with values in the space $ B $ and converging  to zero in this space:

$$
|| \ y_n \ ||B \to 0, \ n \to \infty,
$$
 with the norm

$$
||y||c_0(B) \stackrel{def}{=} \sup_n ||y_n||B.
$$

 \ The equality (2.4) may be reformulated as follows: under what necessary and sufficient  conditions (criterion)

$$
 {\bf P} \left( \{ \xi_n \} \in c_0(B)  \right) = 1. \eqno(2.5)
$$

 \ On the other words,  we raise the question of finding of necessary and sufficient  conditions for
convergence of random elements in (separable) Banach spaces, in the terms of a famous monograph of V.V.Buldygin
 \cite{Buldygin1}. \par

\vspace{4mm}

 \ Let us return to the initial notations: $ \{  f_n(x) \}, \ \nu(\cdot)  $ and so one. Suppose $   f(x) = f_{\infty}(x) = 0. $ \par

 \ We need to introduce some new notations.

$$
\tilde{f}_n(x) :=  \arctan{f_n(x)}, \ \kappa_n^m = \kappa_n^m(F) \stackrel{def}{=}
\int_X \arctan ( \max_{k=n}^m |f_k(x)|) \ \nu(dx) =
$$

$$
\int_X  \max_{k=n}^m |\tilde{f}_k(x)| \ \nu(dx), \ m \ge n + 1; \eqno(2.6)
$$

$$
\overline{\kappa}(F) \stackrel{def}{=} \overline{\lim}_{n \to \infty} \sup_{m \ge n+1}  \kappa_n^m(F). \eqno(2.7)
$$

 \vspace{3mm}

{\bf Definition 2.1.} The introduced above functional $ \overline{\kappa}(F) $
and the like functionals that will be appear further, are sub-linear and
will be named as "critical functional". \par

\vspace{4mm}

{\bf Theorem 2.1.} {\it In order to  }  $  \nu\{x: \lim_{n \to \infty} f_n(x) = 0 \} =1, $ {\it is necessary and sufficient that
 the family $  F  $ belongs to the kernel of the critical functional } $  \overline{\kappa}(F): $

$$
\lim_{n \to \infty} \sup_{m > n} \kappa_n^m (F) = 0 \eqno(2.8)
$$
{\it or equally }

$$
\overline{\kappa}(F) = 0 \hspace{5mm} {\bf or \hspace{4mm} equally} \hspace{5mm} F \in \ker(\overline{\kappa} ).  \eqno(2.8a)
$$

\vspace{3mm}

{\bf Proof. Necessity.}\\

\vspace{2mm}

 \ Let $  \nu \{x: (\lim_{n \to \infty} |f_n(x)| = 0) \} =1,  $  then with at the same value of the $  \nu \ -  $ measure

 $$
 \lim_{n \to \infty} \sup_{m > n} | f_m(x)| = 0
 $$
 and  a fortiori

 $$
 \lim_{n \to \infty} \sup_{m > n} \arctan |f_m(x)| = 0
 $$
$ \nu(\cdot) $  almost surely. We conclude on the basis  of the dominated convergence theorem

$$
\lim_{n \to \infty} \int_X  \sup_{m > n} \arctan |f_m(x)| \ \nu(dx) = 0,
$$
which is quite equivalent to the equality (2.8.) \\

\vspace{4mm}

{\bf Proof. Sufficiency.}\\

\vspace{3mm}

{\bf 0.}  We set  primarily

$$
A = \{ x, x \in X: \lim_{n \to \infty} f_n(x) = 0 \} = \{ x \in X: \lim_{n \to \infty} \tilde{f}_n(x) = 0   \},
$$

 $$
 A_Q = \{ x, x \in X: \forall s = 1,2,\ldots \ \exists N = 1,2,\ldots:  \ \max_{ k \in [N, N + Q]} |\tilde{f}_k(x)| < 1/s \}.
 $$
 Obviously,

 $$
 {\nu}(A) = \lim_{Q \to \infty} {\nu}(A_Q).
 $$

\vspace{3mm}

{\bf 1.} Note first of all that

$$
A = \{ x, \ x \in X: \lim_{n \to \infty} |f_n(x)| = 0   \} = \{ x: \lim_{n \to \infty}  |\tilde{f}_n(x)| = 0 \} =
$$

$$
\cap_s \cup_N  \{ x: \sup_{n \ge N} | \tilde{f}_n(x) |  < 1/s \} \eqno(2.9)
$$
 and correspondingly

 $$
 A_Q := \cap_s \cup_N  \{  x: \max_{n \in [ N, N + Q]} | \tilde{f}_n(x) |  < 1/s \}.
 $$

\vspace{3mm}

{\bf 2.} Let the condition (2.8) be satisfied.  We consider a supplementary  set

$$
B := \overline{A} = X \setminus A, \ B_Q = \overline{A}_Q = X \setminus A_Q. \eqno(2.10)
$$

 \ Elucidation: the set of elementary  events $ B $ may contains (theoretically) also the points when
the limit  does not exists. \par

\ We can write

 $$
 B =  \cup_s \cap_N \{x: \sup_{ n \ge N} | \tilde{f}_n(x)| \ge 1/s  \},\eqno(2.11)
 $$

 $$
 B_Q =  \cup_s \cap_N \{x: \max_{ n \in [N, N+Q]} | \tilde{f}_n(x)| \ge 1/s  \} =
 \cup_s C_{s,Q}, \eqno(2.12)
 $$
where

$$
C_{s,Q} = \cap_N \{x: \max_{ n \in [N, N+Q]} | \tilde{f}_n(x)| \ge 1/s  \} =  \cap_N D^{(N)}_{s,Q}, \eqno(2.13)
$$

$$
D^{(N)}_{s,Q} := \{x, \ x \in X: \ \max_{ n \in [N, N+Q]} | \tilde{f}_n(x)| \ge 1/s  \}. \eqno(2.14)
$$

\vspace{3mm}

{\bf 3.} We obtain using the Tchebychev's  inequality:

$$
 \nu \{ \ \left(  D^{(N)}_{s,Q}  \right) \} \le   \frac{\kappa_N^{N + Q}}{\arctan(1/s)} \to 0, \ N \to \infty,
$$
therefore for all the natural values $ s, Q $

$$
\nu \left( C_{s,Q} \right) =  0,
$$
following

$$
\forall Q = 1,2,\ldots \ \Rightarrow \nu \left(  B_Q  \right) = 0. \eqno(2.15)
$$

\vspace{3mm}

{\bf 4.} We find

$$
 \nu (B) = \lim_{Q \to \infty} \nu \left(B_Q \right)  = 0,\eqno(2.16)
$$
and ultimately $ \nu (A) = 1, $ Q.E.D. \par

 \vspace{3mm}

 {\bf Remark 2.1. } \par

\vspace{3mm}

 \  We will consider  here the case of the Banach space $ c $ consisting on all the numerical sequences $ \{ x(n) \}  $
with existing the limit

$$
\exists \lim_{n \to \infty} x(n) =: x(\infty)
$$
 with at the same norm as above. As before, we consider the classical problem:
  let $ \xi =  \{\xi(n) \} $  be a random sequence; find the conditions (necessary conditions and sufficient
conditions) under which  $ {\bf P} ( \exists \lim_{n \to \infty} \xi(n) ) = 1  $ or equally
$ {\bf P}  ( \{ \xi  \} \in c ) = 1.  $\par
 Notations:

$$
\tilde{\xi}(n) :=  \arctan{\xi(n)}, \ \gamma_n^m = \gamma_n^m(\xi) \stackrel{def}{=}
{\bf E} \arctan ( \max_{k=n}^m |\xi(k) - \xi(n)|).
$$

$$
\overline{\gamma}(\xi) := \overline{\lim}_{n \to \infty} \sup_{m \ge n + 1} \gamma_n^m(\xi).
$$

\vspace{3mm}

{\bf Proposition 2.1.} \  We find analogously to the theorem 2.1:  $ {\bf P} ( \{\xi \} \in c ) =1 $ if and only if

 $$
 \lim_{n \to \infty} \sup_{m \ge n +1} \gamma_n^m (\xi) = 0, \eqno(2.17)
 $$
or briefly

$$
\overline{\gamma}(\xi) = 0, \hspace{5mm} \ {\bf or \hspace{4mm} equally} \hspace{5mm} \xi(\cdot) \in \ker \overline{\gamma}. \eqno(2.17a)
$$

\vspace{3mm}

 {\bf Remark 2.2. }\par

\vspace{3mm}

 \ It is known in the theory of martingales, see e.g. \cite{Burkholder2},  \cite{Burkholder3}, \cite{Hall1},
 chapters 2,3, \cite{Peshkir1} that the estimation of the maximum distribution play a very important role
 for the investigation  of limit theorems, non asymptotical estimations etc. \par
  It follows from our considerations that at the same is true in  more general case of non-martingale
 processes and sequences.

\vspace{3mm}

 {\bf Remark 2.3. } Roughly speaking, the result of theorem 2.1 may be reformulated as follows. Let again
 $ \{ \xi(n) \} $ be a sequence of a r.v., $ \tilde{\xi}(n) = \arctan(|\xi(n)|),  $ and

 $$
 \eta(n) = \sup_{m \ge n} | \tilde{\xi}(m)|.
 $$
  Then

  $$
  \{\omega: \xi(n) \to 0  \}  =  \{\omega: \eta(n) \to 0 \}.
  $$

 \ But the random sequence $ \{ \eta(n) \} $ is monotonically non-increasing,  therefore the sequence
 $ \{ \eta(n) \}  $ tends to zero {\it with probability one} iff this sequence tends to zero {\it in probability,}
 or equally

 $$
 \lim_{n \to \infty} {\bf E} \eta(n) = 0, \eqno(2.18)
 $$
because the variables $ \eta(n) $ are uniformly bounded. \par

\vspace{3mm}

 {\bf Remark 2.4. } Let's turn our attention to the properties of introduced above critical functional
$ \overline{\kappa}(F) $ and $ \overline{\gamma}(\xi) $ and so one. Naturally and obviously, the kernels of these functionals
are {\it closed} linear subspaces;   if for definiteness  $ \overline{\gamma}(\xi_1) = 0 $ and
 $ \overline{\gamma}(\xi_2) = 0, $ then

$$
\overline{\gamma}(c_1 \xi_1 + c_2 \xi_2) = 0, \hspace{4mm} c_1, c_2 = \const. \eqno(2.19)
$$

\vspace{4mm}

\section{Almost everywhere convergence of Fourier series.}

 \vspace{4mm}

   \ Let $  T $  be a segment $ T = [0, 2 \pi] $ and $ \mu   $ is customary {\it renormed } Lebesgue measure
 $  d \mu = dx/(2 \pi). $  \par

 \ Let also $ g = g(x), \ x \in R  $ be measurable and integrable $ (2 \pi) \ - \ $ periodical numerical function.
Denote by $ s_n(x) = s_n[g](x) $ the its $  n^{th}  $  partial sum of ordinary Fourier (trigonometrical) series, which may be
written through  $ (2 \pi) \ - \ $ periodical  convolution  with Dirichlet kernel  $ D_n(x): $

$$
s_n[g](x) = g*D_n(x) = \int_0^{2 \pi} D_n(x-y) \ g(y) \ dy, \eqno(3.1)
$$
where

$$
D_n(x) = \frac{ \sin[(n + 1/2)x]}{ 2 \pi \ \sin(x/2)}. \eqno(3.2)
$$

 \ The problem of finding (sufficient) conditions under the function $  g(\cdot) $ for the almost everywhere
convergence

$$
s_n[g](x) \stackrel{\mu  \ - \ {\bf a.e.}}{\to} g(x)  \eqno(3.3)
$$
 is named ordinary as {\it Luzin's  problem } and has a long history, see for example \cite{Antonov1},
\cite{Carleson1}, \cite{Grafakos1}, \cite{Hunt1}, \cite{Lukomskii1}, \cite{Reyna1}, \cite{Weisz1}. \par

 \ The following famous result belongs to N.Y.Antonov \  \cite{Antonov1}: if

$$
\int_T | \ f(x) \ | \cdot \ \ln^+ | \ f(x) \ |  \cdot  \ln^+  \ln^+  \ln^+ | \ f(x) \ |   \ dx < \infty, \eqno(3.4)
$$
where

$$
\ln^+z = \max(e, \ln z), \hspace{3mm} z > 0,
$$
or equally $  f(\cdot) \in L \ln L \ln \ln \ln L, $ then the convergence (3.3) holds true. \par

\vspace{4mm}

 \ Define also a difference Dirichlet kernel  $ D_{m,n}(x) = D_m(x) - D_{n1}(x), \ m \ge n+1, \ n \ge 2, $ then

$$
D_{m,n}(x) = \pi^{-1} \frac{\sin[(m-n)x/2] \cdot \cos[ (m + n + 1)x/2 ]}{\sin(x/2)},
$$
and introduce the difference of the Fourier   sums

$$
s_{m,n}(x) := [D_{m,n}*g](x) = s_m[g](x) - s_{n}[g](x).
$$
 Denote also

$$
\theta_n^m[g] \stackrel{def}{=} \int_{[0, 2 \pi]} \max_{k \in [n+1,m]}  \arctan   \max_{x \in [0,2 \pi]} |s_{k,n}(x)| \ dx, \eqno(3.5)
$$

$$
\overline{\theta}[g] \stackrel{def}{=} \overline{\lim}_{n \to \infty} \sup_{m \ge n+1} \theta_n^m[g].\eqno(3.5a)
$$

 \ The functional $  g \to \overline{\theta}[g] $ is now the critical sub-linear functional for considered here problem.\par

\vspace{3mm}

 \ It follows immediately from Theorem 2.1 and Proposition 2.1 \\

\vspace{3mm}

{\bf Proposition 3.1.} \  The Fourier sums $  s_n[g](x)  $  for the function $  g = g(x) $ converges almost surely
to the function  $  g = g(x) $ (3.3) if and only if

 $$
 \lim_{n \to \infty} \sup_{m \ge n +1} \theta_n^m [g] = 0. \eqno(3.6)
 $$

{\it or for brevity }

$$
\overline{\theta}(g) = 0 \hspace{5mm} {\bf or \hspace{4mm} equally} \hspace{5mm} g \in \ker(\overline{\theta} ).  \eqno(3.6a)
$$

\vspace{4mm}

\section{ About trial function.}

 \vspace{4mm}

 \ We used in the second section a trial function $  x \to \arctan(|x|).  $ Evidently, it can apply some another functions. \par

\vspace{3mm}

{\bf Definition 4.1.} We will denote by $ KB $  the class of all numerical functions $  \{ \phi \}, \ \phi: R \to R_+  $
satisfying the following conditions:

$$
{\bf A.} \ \phi(x)  \ge 0; \ \phi(x) = 0 \ \Leftrightarrow x = 0,
$$
the condition of positivity;

$$
{\bf B.}  \ 0 < x < y \ \Rightarrow \phi(x) < \phi(y),
$$
the strong monotonicity on the right - hand real axis;\\

\vspace{4mm}

{\bf C.}  Continuity: function $ x \to \phi(x)   $ is continuous on the whole axis $  R. $ \\

\vspace{4mm}

{\bf D.} The function $ \phi(x)  $ is even: $ \forall x \in R  \ \Rightarrow \phi(-x) = \phi(x).   $ \\

\vspace{4mm}

$$
{\bf E.} \ \sup_{x \in R} \phi(x) < \infty,
$$
 the condition of boundedness.\\

 For example:

$$
\phi(x) = \arctan |x|, \ \phi(x) = \frac{|x|}{1 + |x|}, \ \phi(x) = \frac{x^2}{1 + x^2}
$$
and so one. \par

\vspace{3mm}

 \ The  assertion of theorem 2.1 may be rewritten as follows.
Denote as before for any function $  \phi  $ from the set $  KB  $

$$
 \kappa_n^m(\phi) = \kappa_n^m(\phi,F) \stackrel{def}{=}
\int_X \phi ( \max_{k=n}^m |f_k(x)|) \ \nu(dx), \ m \ge n + 1, \eqno(4.1)
$$

$$
\overline{\kappa}(\phi) =
\overline{\kappa}(\phi,F) \stackrel{def}{=} \overline{\lim}_{n \to \infty} \sup_{m \ge n+1} \kappa_n^m (\phi, F).
\eqno(4.1a)
$$

\vspace{4mm}

{\bf Theorem 4.1.}  In order to  $  \nu\{x: \lim_{n \to \infty} f_n(x) = 0 \} =1, $ it is  sufficient that
for {\it some} function  $ \phi(\cdot)  $  from the class  $  KB $

$$
\lim_{n \to \infty} \sup_{m \ge n+1} \kappa_n^m (\phi, F) = 0, \eqno(4.2)
$$
or equally

$$
 \{F \} \in \ker ( \overline{\kappa}(\phi)). \eqno(4.2a)
$$
and is necessary that for {\it arbitrary} function $  \phi  $  belonging to at the same set $  KB  $
the relation (4.2) or (4.2a) holds true. \par

\vspace{3mm}

{\bf Definition 4.2.} We will denote by $ K $  the class of all numerical functions $  \{ \phi \}, \ \phi: R \to R_+  $
satisfying the foregoing conditions at this section {\bf A}, {\bf B}, {\bf C}, {\bf D}, i.e. all the conditions
except for the latter condition of boundedness {\bf E.} For instance, let $  \phi_p(x) := |x|^p, \ p = \const > 0; $ then
$ \phi_p(\cdot) \in K. $\par

 \ In particular, arbitrary Young-Orlicz function  $ \phi(x)  $ belongs to the set $  K. $

\vspace{4mm}

{\bf Theorem 4.2.}  In order to  $  \nu\{x: \lim_{n \to \infty} f_n(x) = 0 \} =1, $ it is  sufficient that
for  some function  $ \phi(\cdot)  $  from the class  $  K $

$$
\lim_{n \to \infty} \sup_{m \ge n+1} \kappa_n^m (\phi, F) = 0, \eqno(4.3)
$$
and  herewith

$$
\lim_{n \to \infty} \int_X \phi( \ ||f_n(x)||B \ ) \ \nu(dx) = 0; \eqno(4.4a)
$$

$$
\sup_n \int_X \phi( \ || \ f_n(x)||B \ ) \ \nu(dx) \le \sup_n \  \sup_{m \ge n+1} \ \kappa_n^m (\phi, F). \eqno(4.4b)
$$

\vspace{4mm}

 \ The relation (4.4a) implies on the language of the theory of Orlicz's function the so-called moment, or weak convergence
$  f_n \to 0 $ in the Orlicz norm $ |||\ \cdot \ |||L\phi; $  the last equality (4.4b) denotes the uniform boundedness of
the considered sequence $ \{  f_n(\cdot)  \}  $ in this space. \par

 \ If in addition this function  $ \phi = \phi(z)  $ is Young-Orlicz function satisfying the so-called $  \Delta_2 $
condition, then the sequence  of functions  $  f_n(\cdot) $ convergent to zero also in the Orlicz's norm
$ |||\ \cdot \ |||L\phi: $

$$
 \lim_{n \to \infty} |||\ f_n \ |||L\phi = 0. \eqno(4.5)
$$

\vspace{3mm}

 \ Let us show another approach which is closely to the so-called Grand Lebesgue Spaces (GLS), see e.g.
 \cite{Fiorenza3}, \cite{Iwaniec2}, \cite{Kozachenko1}, \cite{Liflyand1}, \cite{Ostrovsky1}, \cite{Ostrovsky302}. \par

 \  Let again $ F = \{ f_n(x)\}, \ x \in X $  be as before in the first section be the sequence of measurable functions.
 Define a new function

$$
\psi(p) \stackrel{def}{=} \sup_n \left[  \int_X |f_n(x)|^p \ \nu(dx) \right]^{1/p}, \eqno(4.6)
$$
the so-called {\it  natural } function for the sequence  $ F = \{ f_n(x)\}, \ x \in X, $ and suppose its finiteness
for certain  interval of the form $  1 \le p < R, $ where $ 1 < R = \const \le \infty.  $ \par

 \ This function  $ \psi = \psi(p) $ generated the so - called Grand Lebesgue Space (GLS) $  G\psi  $ consisting
on all the numerical measurable functions $ h = h(x), \ x \in X  $ with finite norm

$$
|| h(\cdot) ||G\psi := \sup_{p \in (1,R) } \left[ \frac{|h|_p}{\psi(p)}  \right],  \eqno(4.7)
$$
where as usually

$$
|h|_p := \left[ \int_X |h(x)|^p \ \nu(dx)   \right]^{1/p}.
$$

 \ These spaces are complete rearrangement invariant Banach spaces which are detail investigated in
\cite{Kozachenko1}, \cite{Liflyand1}, \cite{Ostrovsky1} and so one. \par

\vspace{3mm}

 \ Define the following critical functions

$$
\lambda_n^m(F) = || \ \max_{k = n + 1}^m |f_k(\cdot)| \  ||G\psi,
$$

$$
\overline{\lambda}(F) = \overline{\lim}_{n \to \infty} \sup_{m \ge n+1} \lambda_n^m(F). \eqno(4.8)
$$

\vspace{3mm}

 We conclude as before: \\

 \vspace{3mm}

{\bf Theorem 4.2. }   If $ \overline{\lambda}(F)  = 0,  $ then the sequence $ f_n(x)   $ converges to zero
almost everywhere and in the Grand Lebesgue Space norm $  G \psi. $ \par

\vspace{3mm}

 \ As long as the classical Lebesgue-Riesz spaces $  L_p(X), \ p = \const \ge 1 $
 are the extremal case for GLS spaces, and also the particular cases of
the classical Orlicz spaces,  we conclude denoting again for $  f = \{  f_n(x) \}, \ x \in X   $

$$
\overline{l_p}(F) = \overline{\lim}_{n \to \infty} \sup_{m \ge n+1}  | f_n  - f_m|_p. \eqno(4.9)
$$

\vspace{4mm}

{\bf Theorem 4.3. }   If $ \overline{l_p}(F)  = 0,  $ then the sequence $ f_n(x)   $ converges
almost everywhere as well as in the  Lebesgue-Riesz norm $  L_p(X). $ \par

\vspace{4mm}

\section{ Convergence of random elements in separable Banach spaces.}

 \vspace{4mm}

 \ We return to the raised before the question of finding of necessary and sufficient  conditions for
convergence of random elements in separable Banach spaces, see  \cite{Buldygin1}. \par

 The particular case  of this problem is the problem of continuity of random processes and fields, is considered
in the articles  and books \cite{Dudley1}, \cite{Fernique1}, \cite{Fernique2}, \cite{Fernique3},
 \cite{Garsia1}, \cite{Gikhman1}, \cite{Ibragimov1},  \cite{Hu1},
\cite{Kolmogorov1}, \cite{Kozachenko1}, \cite{Kwapien1}, \cite{Ledoux1}, \cite{Nisio1},
\cite{Ostrovsky1}, \cite{Ostrovsky301}, \cite{Slutsky1}, \cite{Sudakov1},
\cite{Talagrand1}, \cite{Talagrand2}, \cite{Vinkler1}, \cite{Wakhaniya1},
\cite{Watanabe1} etc.\par

 \ This problem may be easily reduced to the problem of uniform convergence of random numerical functions, for
instance, the problem of uniform convergence of expression of the series by Franklin  orthogonal  system,
see \cite{Ciesielski1}, \cite{Franklin1},  \cite{Ostrovsky301}. \par

\ We consider here the problem of convergence of random elements with values in Banach space.\par

 \ In detail: let $ \zeta =  \{\zeta(n) \} $  be a random sequence with values in the separable Banach space $  B; $
find the conditions (necessary conditions and sufficient
conditions) under which  $ {\bf P} ( \exists \lim_{n \to \infty} \zeta(n) ) = 1  $ or equally
$ {\bf P}  ( \{ \zeta  \} \in c(B) ) = 1.  $\par
 \ Notations:

$$
\tilde{\zeta}(n) :=  \arctan{\zeta(n)}, \ \tau_n^m = \tau_n^m(\xi) \stackrel{def}{=}
{\bf E} \arctan ( \max_{k=n}^m || \ \zeta(k) - \zeta(n) \ ||B), \eqno(5.1)
$$

$$
\overline{\tau}(\zeta) := \overline{\lim}_{n \to \infty} \sup_{m \ge n + 1} \tau_n^m(\zeta). \eqno(5.1a)
$$

\vspace{3mm}

{\bf Proposition 5.1.} \  We conclude analogously to the theorem 2.1:  $ {\bf P} ( \{\zeta \} \in c ) =1 $ if and only if

 $$
 \lim_{n \to \infty} \sup_{m \ge n +1} \tau_n^m (\xi) = 0, \eqno(5.2)
 $$
or briefly

$$
\overline{\tau}(\zeta) = 0, \hspace{5mm} \ {\bf or \hspace{4mm} equally} \hspace{5mm} \zeta(\cdot) \in \ker \overline{\tau}. \eqno(5.3)
$$

\vspace{4mm}

\section{Concluding remarks.}

\vspace{4mm}

 \hspace{2mm} {\bf I.} The case of  metric (linear) space  $  B. $ \par

 \hspace{1mm} It is not hard  to generalize obtained above results on the case when the Banach space $  B  $ is replaced by certain
separable linear metric space $ L $ equipped with translation invariant metric function  $ \rho = \rho(x - y).   $\par

\vspace{4mm}

{\bf II.} Recall that for the convergence in probability (measure) of the sequence of Banach space valued r.v. $  \eta_n  $
 the necessary and sufficient condition is following

$$
\lim_{n,m \to \infty} {\bf E} \arctan || \ \eta_n - \eta_m \ ||B = 0.
$$

\vspace{4mm}

{\bf III.} It is interest by our opinion to investigate in the spirit of this article
the case of the non-sequential  convergence; as well as to obtain the criterion for a.e.
convergence for multiple sequences, especially multiple Fourier series.\\

\end{document}